\documentclass[12pt]{amsart}

\usepackage{amssymb}
\usepackage{mathrsfs}
\usepackage{amscd}
\usepackage[numbers]{natbib}
\usepackage{amsfonts}

\usepackage{color}
\usepackage[T1]{fontenc}

\usepackage[scaled=.95]{helvet}
\usepackage{courier}

\usepackage[english]{babel}
\usepackage[utf8]{inputenc}
\usepackage{latexsym}   
\usepackage{tikz}
\usepackage{amsthm,times}
\usepackage[margin=1in]{geometry}
\usepackage{epsfig,graphicx}
\usepackage[matrix,arrow]{xy}
\usepackage{subfigure}
\usepackage{mathrsfs}
\usepackage{mathtools}
\usepackage{graphicx,color}
\usepackage{eulervm}
\usepackage{setspace}

\newcommand{\lar}{\longrightarrow}
\newcommand{\llar}{-\kern-5pt-\kern-5pt\longrightarrow}

\newtheorem{Theorem}{Theorem}[section]
\newtheorem{Lemma}[Theorem]{Lemma}

\newtheorem{Proposition}[Theorem]{Proposition}

\newtheorem{Example}[Theorem]{Example}
\newtheorem{Conjecture}[Theorem]{Conjecture}
\newtheorem{Definition}[Theorem]{Definition}

\def\sqr#1#2{{\vcenter{\hrule height.#2pt
			\hbox{\vrule width.#2pt height#1pt \kern#1pt
				\vrule width.#2pt}
			\hrule height.#2pt}}}
\def\phi{\varphi}

\DeclareMathOperator{\Hom}{Hom}

\DeclareMathOperator{\Image}{Im}

\DeclareMathOperator{\coker}{coker}

\DeclareMathOperator{\Ann}{Ann}

\DeclareMathOperator{\Supp}{Supp}
\DeclareMathOperator{\Ass}{Ass}

\DeclareMathOperator{\Min}{Min}

\DeclareMathOperator{\ext}{Ext}

\DeclareMathOperator{\Tor}{Tor}

\DeclareMathOperator{\Ht}{ht}

\def\Ree#1{{\cal R}(#1)}

\def\ker{{\rm ker}\,}

\def\restr{{\kern-1pt\restriction\kern-1pt}}

\def\R{{\mathbb R}}
\def\C{{\mathbb C}}

\begin{document}
	\newcommand{\mmbox}[1]{\mbox{${#1}$}}
	\newcommand{\proj}[1]{\mmbox{{\mathbb P}^{#1}}}
	\newcommand{\Cr}{C^r(\Delta)}
	\newcommand{\CR}{C^r(\hat\Delta)}
	\newcommand{\affine}[1]{\mmbox{{\mathbb A}^{#1}}}
	\newcommand{\caps}[3]{\mmbox{{#1}_{#2} \cap \ldots \cap {#1}_{#3}}}
	\newcommand{\N}{{\mathbb N}}
	\newcommand{\im}{\mathop{\rm Im}\nolimits}
	\newcommand{\arrow}[1]{\stackrel{#1}{\longrightarrow}}
	\newcommand{\CB}{Cayley-Bacharach}
	\newtheorem{defn0}{Definition}[section]
	\newtheorem{prop0}[defn0]{Proposition}
	\newtheorem{quest0}[defn0]{Question}
	\newtheorem{thm0}[defn0]{Theorem}
	\newtheorem{lem0}[defn0]{Lemma}
	\newtheorem{corollary0}[defn0]{Corollary}
	\newtheorem{example0}[defn0]{Example}
	\newtheorem{remark0}[defn0]{Remark}
	\newtheorem{prob0}[defn0]{Problem}
	
	\newenvironment{defn}{\begin{defn0}}{\end{defn0}}
	\newenvironment{prop}{\begin{prop0}}{\end{prop0}}
	\newenvironment{quest}{\begin{quest0}}{\end{quest0}}
	\newenvironment{thm}{\begin{thm0}}{\end{thm0}}
	\newenvironment{lem}{\begin{lem0}}{\end{lem0}}
	\newenvironment{cor}{\begin{corollary0}}{\end{corollary0}}
	\newenvironment{exm}{\begin{example0}\rm}{\end{example0}}
	\newenvironment{rem}{\begin{remark0}\rm}{\end{remark0}}
	\newenvironment{prob}{\begin{prob0}\rm}{\end{prob0}}
	
	\newcommand{\defref}[1]{Definition~\ref{#1}}
	\newcommand{\propref}[1]{Proposition~\ref{#1}}
	\newcommand{\thmref}[1]{Theorem~\ref{#1}}
	\newcommand{\lemref}[1]{Lemma~\ref{#1}}
	\newcommand{\corref}[1]{Corollary~\ref{#1}}
	\newcommand{\exref}[1]{Example~\ref{#1}}
	\newcommand{\secref}[1]{Section~\ref{#1}}
	\newcommand{\remref}[1]{Remark~\ref{#1}}
	\newcommand{\questref}[1]{Question~\ref{#1}}
	\newcommand{\probref}[1]{Problem~\ref{#1}}

	\newcommand{\std}{Gr\"{o}bner}
	\newcommand{\jq}{J_{Q}}
	\def\Ree#1{{\mathcal R}(#1)}
	
	



\title[On the factorial case of Huneke's conjecture]{On the factorial case of Huneke's conjecture for local cohomology modules}

\author[Dosea, \, Miranda-Neto]{Andr\'e Dosea, \, Cleto B. Miranda-Neto}

\address{Departamento de Matem\'atica, Universidade Federal da Para\'iba - 58051-900, Jo\~ao Pessoa, PB, Brazil}
\email{andredosea@hotmail.com}
\address{Departamento de Matem\'atica, Universidade Federal da Para\'iba - 58051-900, Jo\~ao Pessoa, PB, Brazil}
\email{cleto@mat.ufpb.br}



\begin{abstract} A conjecture raised in 1990 by C. Huneke predicts that, for a $d$-dimensional Noetherian local ring $R$, local cohomology modules of finitely generated $R$-modules have finitely many associated primes. Although counterexamples do exist, the conjecture has been confirmed in several cases, for instance if $d\leq 3$, and witnessed some progress in special cases for higher $d$. In this paper, assuming that $R$ is a factorial domain, we establish the case $d=4$, and under different additional conditions (in a couple of results) also the case $d=5$. Finally, when $R$ is regular and contains a field, we apply the Hartshorne-Lichtenbaum vanishing theorem as a tool to deal with the case $d=6$.
\end{abstract}



\maketitle
\section{Motivation: Huneke's conjecture}

In this paper we are motivated by the following conjecture proposed by Craig Huneke in 1990. 

\begin{Conjecture}{\rm (\cite[Conjecture 5.1]{Hu})} Let $R$ be a Noetherian local ring and $M$ a finitely generated $R$-module. Then, for each $i\geq 0$ and any ideal $I$ of $R$, the $R$-module ${\rm H}^i_I(M)$ has finitely many associated primes.
    
\end{Conjecture}

Here, $H^i_I(M)=\varinjlim {\rm Ext}_R^i(R/I^n,M)$, the $i$-th local cohomology module of $M$ with respect to the ideal $I$, which is a classical object in commutative algebra and algebraic geometry (see, e.g., \cite{Brodman}). It is well-known that if the set of associated primes of a given local cohomology module is finite then its support is Zariski-closed.

The conjecture has attracted the attention of numerous researchers and has been
answered positively in several cases, typically in small dimensions; see \cite{Bah2}, \cite{Hu-Koh}, \cite{Marley1}, \cite{Mel}. We also refer to \cite{HuS} and \cite{Lyub}. It is worth recalling that there are counterexamples to Huneke's conjecture (see \cite{Moty}, \cite{SS}). Additionally, if $R$ is not local, these sets of associated
primes may be infinite (see \cite{Sin}), whereas they must be finite if $R$ is a
smooth ${\mathbb Z}$-algebra (see \cite{multi}). 

In the present paper, we are first concerned with the cases where ${\rm dim}\,R$ is either 4 or 5, under certain additional assumptions which include $R$ being  factorial (i.e., a unique factorization domain). Our result in the former case is Theorem \ref{four dimensional case}, while in the latter our main results are Theorem \ref{five dimensional hartshorne} and Theorem \ref{five dimensional part huneke}. Finally, if $R$ is a regular local ring containing a field, we treat in Theorem  \ref{theorem dim 6} the case ${\rm dim}\,R=6$; to the best of our knowledge, this is the first existing result towards the six-dimensional case of Huneke's conjecture. We provide explicit examples in all three situations.

Throughout this note, all rings are assumed to be commutative and Noetherian, and by a {\it finite} module we mean a finitely generated module. Wherever they appear, $R$ denotes a ring and $I$ a proper ideal of $R$, whose height we write $\Ht I$. As usual, we set $V(I)=\{P\in {\rm Spec}\,R \mid I\subset P\}$.

\section{Auxiliary notions and facts}\label{section preliminarues}

We present some concepts that are central in this note. Details can be found in \cite{Mafi1}.

	\begin{Definition}\label{weakly}\rm
Let $H$ be an $R$-module.
\begin{itemize}
    \item[\rm (i)] $H$ is \textit{weakly Laskerian} if the set $\Ass_R H/U$ is finite for every $R$-submodule $U$ of $H$.
    \item[\rm (ii)] $H$ is {\it $I$-weakly cofinite} if $\Supp_RH \subset V(I)$ and if the $R$-module $\ext^i_R (R/I, H)$ is weakly Laskerian for all $i \geq 0$.
\end{itemize} 
\end{Definition}

The next example and Proposition 
\ref{basic properties weakly laskerian} collect useful properties which we shall freely use  without explicit mention.

\begin{Example}\rm \label{weakly laskerian examples}
All Noetherian modules are weakly Laskerian.
If $\Supp_R H$ is finite, then $H$ is weakly Laskerian; consequently, all Artinian modules are weakly Laskerian. If ${\rm Supp}_RH \subset V(I)$ and $H$ is weakly Laskerian, then $H$ is $I$-weakly cofinite.
\end{Example}



\begin{Proposition}\label{basic properties weakly laskerian}
The following assertions hold:
\begin{itemize}
    \item[\rm (i)] If $0 \to N \to M \to L \to 0$ is a short exact sequence of $R$-modules, then $M$ is weakly Laskerian $($resp. $I$-weakly cofinite$)$ if and only if $N$ and $L$ are weakly Laskerian $($resp. $I$-weakly cofinite$)$;

    \item[\rm (ii)] If $M$ is finite and $N$ is weakly Laskerian, then $\Tor^R_i (M,N)$ and $\ext^i_R (M,N)$ are weakly Laskerian for all $i \geq 0$;
    \item[\rm (iii)] If $M$ is $I$-weakly cofinite, then $\Ass_R M$ is finite. 
    \end{itemize}
    \end{Proposition}
    
    
In order to address Huneke's conjecture, our main strategy in this paper (with the exception of the results in Section \ref{sixdim}, where another approach is used) is to investigate when local cohomology modules with support in $I$ are $I$-weakly cofinite, and then to apply property (iii) above.

Note that Definition \ref{weakly}(ii) extends the classical notion of $I$-cofiniteness, i.e., $H$ is {\it $I$-cofinite} if $\Supp_RH \subset V(I)$ and $\ext^i_R (R/I, H)$ is a finite $R$-module for all $i \geq 0$.  The lemma below, which will be used later on, provides a source of examples of $I$-cofinite local cohomology modules.

    \begin{Lemma}{\rm (\cite[Theorem 1]{Kawazaki})}\label{lemma principal ideal}
Let $M$ be a finite $R$-module and $I$ a principal ideal of $R$. Then, $H^i_I(M)$ is $I$-cofinite for all $i\geq 0$.    
\end{Lemma}

In particular, under the conditions of Lemma \ref{lemma principal ideal}, the $R$-module $H^i_I(M)$ must have finitely many associated primes for each $i\geq 0$.   

Finally, we recall a classical result on local cohomology that we will apply throughout without explicit mention, namely, Grothendieck's vanishing theorem (see, e.g., \cite[Theorem 6.1.2]{Brodman}), which asserts that if $M$ is an $R$-module then $$H_I^i(M)=0 \quad \mbox{for\, all} \quad i> {\rm dim}\,M$$ (note that if $M$ is finite then ${\rm dim}\,M$ is just the dimension of the quotient ring $R/{\rm ann_RM}$). In particular, if $P$ is a prime ideal of $R$ containing $I$, then $H_{I_P}^j(R_P)=0$ for all $j> \Ht P$. Also, it is well-known that if $R$ is local and $M$ is finite then $\Ass_R H^i_I (M)$ is finite for $i=0, 1$ (see, e.g.,
\cite[Proposition 1.1(c)]{Marley1}), so that for Huneke's conjecture the interest lies in the case $i\geq 2$.

\section{The four-dimensional factorial case}

First we consider the following lemma (in arbitrary dimension).

        \begin{Lemma}{\rm (\cite[Theorem 3.1]{Bah2})}\label{lemma dim menor que 2}
         Let $R$ be a local ring and $M$ a finite $R$-module. If $t$ is a non-negative integer such that $\dim H^i_I (M) \leq 2$ for all $i<t$, then $H^i_I (M)$ is $I$-weakly cofinite for all $i<t$ and  $\Hom_R (R/I, H^t_I (M))$ is weakly Laskerian.   
        \end{Lemma}

        Now we are in a position to present the main result of this section. Our theorem substantially improves \cite[Theorem 3.6]{Bah3}, where $R$ was taken  to be a {\it regular} local ring.
        
    \begin{Theorem}\label{four dimensional case} Let $R$ be a four-dimensional factorial local ring and $M$ a finite $R$-module. Then, $H^i_I (M)$ is $I$-weakly cofinite for all $i$. In particular, $\Ass_R H^i_I (M)$ is finite for all $i$.
    \end{Theorem}
    \begin{proof} We can suppose without loss of generality that $I$ is a radical ideal. Also notice that we may assume $\Ht I \leq 1$. Indeed, the case $\Ht I \geq 2$ can be settled by using Lemma \ref{lemma dim menor que 2} since $\dim R/I \leq 2$ in this case, which implies $\dim H^i_I (M) \leq 2$ for all $i$.

     Now we denote $X = \{P \in \Min_R R/I \mid \Ht P \leq 1\}$ and $Y = \{ Q \in \Min_R R/I \mid \Ht Q \geq 2\}$.
	    If $Y$ is empty, then (since $R$ is factorial) $I$ is necessarily a principal ideal and the assertion follows from
        Lemma \ref{lemma principal ideal}. Therefore, we may assume $Y$ is non-empty. 
        
        Let us consider the ideals
        $$J= \displaystyle\bigcap_{P \in X} P \quad \mbox{and} \quad K = \displaystyle\bigcap_{Q \in Y} Q.$$
         Given a prime  ideal $P'$ containing $J+K$, we have  $P' \supset P$ for some $P \in X$ and also $P' \supset Q$ for some $Q \in Y$. This shows that $\Ht P' \geq 2$. If $\Ht P'=2$, then $P' =Q$ and we conclude that $Q= P$ by minimality, which is a contradiction. So $\Ht P'\geq 3$, which (as ${\rm dim}\,R=4$) forces 
     $$
     V(J+K) = \Min_R R/(J+K) \cup \{\textbf{m}\}
     $$ (where $\textbf{m}$ stands for the maximal ideal of $R$) and therefore $V(J+K)$ 
    is a finite set.
        Now, the Mayer-Vietoris sequence yields the exact sequence
        $$
        H^i_{J+K} (M) \to H^i_J (M) \oplus H^i_K (M) \stackrel{f}\lar H^i_I (M) \to H^{i+1}_{J+K} (M).
        $$
        Note that, as $R$ is a factorial domain, the ideal $J$ must be principal and then each $H^i_J (M)$ is $J$-cofinite by Lemma \ref{lemma principal ideal}. Moreover, since $$\dim R/K \leq 4-\Ht K\leq  2,$$ we can apply  Lemma \ref{lemma dim menor que 2} to obtain that $H^i_K (M)$ is $K$-weakly cofinite.
        Furthermore, because $V(J+K)$ is  finite, the modules $H^i_{J+K} (M)$ and $H^{i+1}_{J+K} (M)$ have finite support and hence are weakly Laskerian.
        The same is valid for the modules $\ext^j_R (R/J, H^i_K (M))$ and $\ext^j_R (R/K, H^i_J (M))$.
        
Next, we consider the exact sequence $$0 \to V \to H^i_J (M) \oplus H^i_K (M) \stackrel{f}\lar H^i_I (M) \to W \to 0,$$ where the $R$-modules $V={\rm Ker}\,f$ and $W={\rm Coker}\,f$  are seen to be weakly Laskerian. It splits into the two short exact sequences
        $$0 \to V \to H^i_J (M) \oplus H^i_K (M) \to \Image f \to 0,$$
        $$0 \to \Image f \to H^i_I (M) \to W \to 0.$$ Applying to these sequences the functors $\ext^j_R (R/J, -)$ and $\ext^j_R (R/K, -)$,
        we derive the weak Laskerianess of the $R$-modules $$\ext^j_R (R/J, H^i_I (M)) \quad \mbox{and} \quad \ext^j_R (R/K, H^i_I (M)).$$
        Finally, applying the functor $\ext^j_R (-, H^i_I (M))$ to the standard exact sequence $$0 \to R/ I \to R/J \oplus R/K \to R/(J+K) \to 0,$$ we deduce the $I$-weak cofiniteness of $H^i_I (M)$.
	\end{proof}

 \begin{Example}\rm
    Let us begin with the regular local ring $S =k[x_1, x_2, x_3, x_4, x_5]_{(x_1, x_2, x_3, x_4, x_5)}$, where $x_1, x_2, x_3, x_4, x_5$ are indeterminates over a field $k$. Write ${\bf m}=(x_1, x_2, x_3, x_4, x_5)S$ and consider $$I = ((x_1)\cap(x_2, x_3))S \quad \mbox{and} \quad P= (x_1, x_2, x_3)S.$$ According to \cite[Example 3.8]{Marley}, the $S_P$-module $\Hom_S (S/I, H^2_I (S))_P$ is not finite, and hence $\Hom_S (S/I, H^2_I (S))_Q$ is not a finite $S_Q$-module for any $Q\in V(P)$. On the other hand, if we choose $Q\neq {\bf m}$, the module  $H^i_I (S)_Q$ is $I_Q$-weakly cofinite for all $i$. Indeed, if $Q=P$ (in which case ${\rm dim}\,R_Q=3$), this assertion follows from \cite[Corollary 5.3]{Goto}, and if $P \subsetneqq Q$ we apply Theorem \ref{four dimensional case} to the four-dimensional regular (hence factorial) local domain $R=S_Q$.
    In particular, $$\Hom_S (S/I, H^2_I (S))_Q$$ is a weakly Laskerian  (though not finite) $R$-module.
    \end{Example}

    \begin{Example}\rm
    Consider the four-dimensional local ring $$R = k[x_1, x_2, x_3, x_4, x_5]_{(x_1, x_2, x_3, x_4, x_5)}/(x_1^2+ x_2^2+x_3^2+x_4^2+x_5^2),$$ where $x_1, x_2, x_3, x_4, x_5$ are indeterminates over a field $k$ such that either $k=\C$ or $k\subset \R$. Write $q=x_1^2+ x_2^2+x_3^2+x_4^2+x_5^2$, which is a non-degenerate quadratic form. Then, the Klein-Nagata Theorem (see \cite[Theorem 8.2]{Samuel}) guarantees that $R$ is a factorial domain. Consider the ideal $$I= ((x_1,x_2) \cap (x_3,x_4,x_5))R \subset (x_3,x_4,x_5)R=(x_1^2+x_2^2, x_3,x_4,x_5)/(q),$$ which satisfies ${\rm dim}\,R/I=2$. Now let $Q$ be the ideal of $R$ given by \[ Q=\left\{ \begin{array}{cr}
    (x_1 \pm ix_2, x_3, x_4, x_5)R, ~~ if ~~ k=\C\\\\
    (x_1^2+x_2^2, x_3,x_4,x_5)R, ~~ if ~~ k\subset \R
    \end{array} \right.
  \] and note that ${\rm dim}\,R/Q=1$ and $Q\in {\rm Min}_RR/I$. Further, it is clear (e.g., by the well-known Jacobian criterion) that the hypersurface local domain $R$ is analytically normal. Such facts put us in a position to apply \cite[Proposition 3.6]{Marley}, which yields that $$\Hom_R (R/I, H^3_I (R))$$ is not a finite $R$-module. On the other hand, this $R$-module is weakly Laskerian by virtue of Theorem \ref{four dimensional case}. It is worth observing that
    this particular example cannot be deduced from \cite[Theorem 3.6]{Bah3} because $R$ is not a regular ring.
    \end{Example}

  \section{Two theorems on the five-dimensional factorial case}\label{twothm}      In the previous section, we applied Lemma \ref{lemma dim menor que 2} as a key tool to settle the four-dimensional factorial case of Huneke's conjecture.  We may wonder whether a similar strategy holds in the five-dimensional case as well. As we will clarify, this can be handled with appropriate (and rather mild) additional assumptions 
  
  
\subsection{First theorem} The first main result of this section is the following.


\begin{Theorem}\label{five dimensional hartshorne}
    Let $R$ be a five-dimensional factorial local ring and $M$ a finite $R$-module whose annihilator is not contained in any minimal prime divisor of $I$ {\rm (}e.g., if ${\rm ann}_RM$ contains an $R/I$-regular element{\rm )}. If either $\Ht I \neq 1$ or $I$ is a principal ideal, then $H^i_I (M)$ is $I$-weakly cofinite for all $i$. In particular, $\Ass H^i_I (M)$ is finite.
\end{Theorem}
\begin{proof} In view of Lemma \ref{lemma principal ideal}, we may assume $\Ht I \geq  2$. Consider first the case  $\Ht I \geq  3$. Then
   $$\dim H^i_I (M) \leq \dim R/I \leq 5 - \Ht I \leq 2 \quad \mbox{for\, all} ~~ i,$$  and so Lemma \ref{lemma dim menor que 2} gives the result. So it remains to investigate the case $\Ht I = 2$. We claim that, also in this situation,  $\dim H^i_I (M) \leq 2$ for all $i$. Indeed, if $P \in \Supp_R H^i_I (M)$ is such that $\Ht P =2$, then 
$$P \in \Min_R R/I \cap \Supp_R M=\Min_R R/I \cap V({\rm ann_RM}),
$$ which is a contradiction.
Hence $\Ht P \geq 3$ for all $P \in \Supp_R H^i_I (M)$. This implies the inequality $\dim H^i_I (N) \leq \dim R/P \leq 2$ for all $i$, as claimed.  Once again, the assertion follows by Lemma \ref{lemma dim menor que 2}. 
\end{proof}

It now seems natural to ask:

\begin{quest} \rm Does Theorem \ref{five dimensional hartshorne} hold true if $I$ is non-principal and $\Ht I=1$?
    
\end{quest}


\subsection{Second theorem} In order to allow height one non-principal ideals in Theorem \ref{five dimensional hartshorne}, we shall restrict ourselves to the situation where $\dim M \leq 4$ in our second main result toward the five-dimensional case (Theorem \ref{five dimensional part huneke}). While, so far, we have investigated the finiteness of $\Ass_R H_I^i (M)$ by means of $I$-weak cofiniteness, the approach to be developed here still makes use of this property but in a slightly different way; a key step is proving that a certain category is Abelian (Proposition \ref{theorem abelian category}), and to this end an ingredient is the following proposition, which is also of self-interest. 

\begin{Proposition}\label{proposition criteria dim < 2}
	Let $R$ be a local ring and $M$ an $R$-module with $\dim M \leq 2$. The following assertions are equivalent:
	\begin{itemize}
	    \item[(\rm i)] $H^i_{I} (M)$ is $I$-weakly cofinite for all $i \geq 0$;
	    \item[(\rm ii)] $\ext^i_R (R/I,M)$ is weakly Laskerian for all $i \geq 0$;
	    \item[(\rm iii)] ${\rm Hom}_R (R/I,M)$ and $\ext^1_R (R/I,M)$ are weakly Laskerian.
	\end{itemize}
	\end{Proposition}
	\begin{proof} Applying \cite[Proposition 3.9]{Mel} for the
	 category of weakly Laskerian $R$-modules, we see that (i) implies (ii).
	It is clear that (ii) implies (iii). 
	Now, assume (iii) and let us prove that (i) holds. 
	
	By flat base change (see \cite[Theorem 4.3.2]{Brodman}) along with \cite[Lemma 2.1]{Marley}, we can assume that the local ring $R$ is ${\bf m}$-adically  complete, where ${\bf m}$ is its maximal ideal. Given two non-negative integers $i,j$, we consider the $R$-module $$T=\ext^i_R (R/I, H^j_I (M)).$$ For any given $R$-submodule $L$ of $T$, we need to show that $\Ass_R T/L$ is finite. Suppose, by way of contradiction, that this set is infinite for some $L$. So, we can choose  a countably infinite subset $$\displaystyle\ \{P_t\}_{t=1}^{\infty}\subset \Ass_R T/L, \quad P_t \neq {\bf m} \quad \mbox{for \,each} ~~ t.$$ Consider the multiplicative closed subset $S = R\setminus \bigcup_{t=1}^{\infty} P_t$. Then, for every $t$, we have $S^{-1}P_t \in \Ass_{S^{-1}R} S^{-1}T/S^{-1}L$. Notice, in particular, that  $\Ass_{S^{-1}R} S^{-1}T/S^{-1}L$ is infinite. On the other hand, \cite[Lemma 3.2]{Marley} (where $R$ is required to be complete) ensures that ${\bf m}\nsubseteq \bigcup_{t=1}^{\infty} P_t$, i.e,  ${\bf m}$ meets $S$, which (as ${\rm dim}\,M\leq 2$) forces $$\dim  S^{-1}M \leq 1.$$ Now, since the $S^{-1}R$-module $\ext^i_{S^{-1}R} (S^{-1}R/S^{-1}I, S^{-1}M)$ is  weakly Laskerian for $i=0,1$, we are in a position to apply \cite[Proposition 3.4]{Roshan}  to deduce that $H^j_{S^{-1}I} (S^{-1}M)$ is $S^{-1}I$-weakly cofinite for all $j \geq 0$. It follows that   $S^{-1}T$ is a weakly Laskerian $S^{-1}R$-module and hence so is the quotient $S^{-1}T/S^{-1}L$. In particular, $\Ass_{S^{-1}R} S^{-1}T/S^{-1}L$ is finite, which is a contradiction.
	\end{proof}

 According to the terminology introduced in \cite{Roshan}, an $R$-module $M$ is said to be {\it $I$-{\rm ETH} weakly cofinite} if $\ext^i_R (R/I, M)$ is weakly Laskerian for all $i\geq 0$.

 \begin{Proposition}\label{theorem abelian category}
    Let $R$ be a local ring. Then, the category of all $I$-{\rm ETH} weakly cofinite modules of dimension at most $2$ is Abelian. In particular, the category of all $I$-weakly cofinite modules of dimension at most $2$ is Abelian.
    \end{Proposition}
    \begin{proof}
    Let $f:M \rightarrow N$ be an $R$-linear map, where $M, N$ are $I$-ETH weakly cofinite modules of dimension at most 2. So, the modules $\ker f$, $\Image f$ and $\coker f$ have dimension at most 2 as well. The monomorphism
    $$\Hom_R (R/I, \ker f) \to \Hom_R (R/I, M)$$ shows that $\Hom_R (R/I, \ker f)$ is weakly Laskerian. Similarly, the module $\Hom_R (R/I, \Image f)$ is weakly Laskerian.
    On the other hand, the short exact sequence $0 \to \ker f \to M \to \Image f \to 0$ yields the exact sequence
    $$\Hom_R (R/I, \Image f) \to \ext^1_R (R/I, \ker f) \to \ext^1_R (R/I, M)$$
    from which we obtain the weak Laskerianess of $\ext^1_R (R/I, \ker f)$.
    Hence, by Proposition \ref{proposition criteria dim < 2}, $\ker f$ is $I$-ETH weakly cofinite.
    Furthermore, the exact sequence
    $$
    \ext^1_R (R/I, M) \to \ext^1_R (R/I, \Image f) \to \ext^2_R (R/I,\ker f)
    $$ allows us to deduce the weak Laskerianess of $\ext^1_R (R/I, \Image f)$. Now, another application of Proposition \ref{proposition criteria dim < 2} shows that $\Image f$ is $I$-ETH weakly cofinite.
    Finally, using the short exact sequence
    $$0 \to \Image f \to N \to \coker f \to 0,$$ we easily deduce that $\coker f$ is $I$-ETH weakly cofinite.
    \end{proof}

It is worth mentioning that, in case $R$ is a five-dimensional unramified {\it regular} local ring and $M$ is a torsion-free finite $R$-module, each set $\Ass H^i_I (M)$ has been proven to be finite (see \cite[Theorem 2.11]{Marley1}). Our Theorem \ref{five dimensional part huneke} below will be established in another direction. Under appropriate assumptions on $M$, which will be no longer required to be torsion-free, we shall deal with the larger class  of factorial local domains (as in Theorem \ref{four dimensional case} and Theorem \ref{five dimensional hartshorne}).

\begin{Theorem}\label{five dimensional part huneke} Let $R$ be a five-dimensional factorial local ring and $M$ a finite $R$-module whose annihilator is not contained in any minimal prime divisor of $I$ {\rm (}e.g., if ${\rm ann}_RM$ contains an $R/I$-regular element{\rm )}. If $\dim M \leq 4$, then $\Ass_R H^i_I (M)$ is finite for all $i$.
     \end{Theorem}
\begin{proof} First, by \cite[Corollary 2.7]{Marley1}, we may assume $\dim M = 4$. Note $H^j_I (M)=0$ for all $j\geq 5$.
By virtue of \cite[Corollary 2.5]{Marley1} (resp.\,\cite[Lemma 2.6]{Cuong}), the $R$-module $H^4_I (M)$ (resp.\,$H^3_I (M)$) has finite support, hence finitely many associated primes. Therefore, in view of Theorem \ref{five dimensional hartshorne}, it remains to prove the finiteness of $\Ass_R H^2_I (M)$ when $\Ht I=1$. To this end, set $X,Y,J$ and $K$ as in the proof of Theorem \ref{four dimensional case}, and similarly to that proof we may suppose $Y$ is non-empty. The Mayer-Vietoris sequence yields the exact sequence
\begin{equation}\label{eq.1}
 H^2_{J+K} (M) \stackrel{f}\lar H^2_J (M) \oplus H^2_K (M) \to H^2_I (M) \to H^3_{J+K} (M) \stackrel{g}\lar H^3_J (M) \oplus H^3_K (M).
\end{equation}
According to \cite[Lemma 2.6]{Cuong}, all modules on the right of $H^2_I (M)$ in (\ref{eq.1}) have finite support. In particular, they are weakly Laskerian, which implies the weak Laskerianess of $\ker g$. This forces $\ker g$ to be $(J+K)$-ETH weakly cofinite.
Our goal now is to show that $\coker f$ has the same property. We argue as follows. First, since $J$ is a principal ideal, we get that $H^i_J (M)$ is $J$-cofinite by Lemma \ref{lemma principal ideal}. Moreover, as $\Ht K \geq 2$ and the intersection $$\Min_R R/K \cap \Supp_R M=\Min_R R/K \cap V({\rm ann_RM})$$ is empty, Theorem \ref{five dimensional hartshorne} yields the $K$-weak cofiniteness of $H^i_K (M)$ for all $i$. On the other hand, since $\Ht (J+K) \geq 3$, we have $$\dim H^i_{J+K}(M) \leq \dim R/(J+K)\leq 2$$ and hence
Lemma \ref{lemma dim menor que 2} allows us to conclude that $H^i_{J+K} (M)$ is $(J+K)$-weakly cofinite. We can apply \cite[Lemma 2.2]{Roshan} to deduce that all modules on the left side of $H^2_I (M)$ in (\ref{eq.1}) are $(J+K)$-ETH weakly cofinite. 
In addition, it follows from \cite[Theorem 1.1(ii)]{Saremi} that all modules in \eqref{eq.1} have dimension bounded above by 2. So, as an application of Proposition \ref{theorem abelian category}, we deduce that $\coker f$ is $(J+K)$-ETH weakly cofinite. 

Now, by means of the short exact sequence
$$0 \to \coker f \to H^2_I (M) \to \ker g \to 0$$ we conclude that
$\Hom_R (R/(J+K), H^2_I (M))$ is weakly Laskerian. Therefore, the set
     $$
     \Ass_R H^2_I (M) \cap V(J+K)
     $$ is finite. Let $P \in \Ass_R H^2_I (M)$. Since $I \subset P$, we obtain $J \subset P$ or $K \subset P$.
     If both inclusions are true, then $P \in V(J+K)$. Otherwise, the exact sequence \eqref{eq.1} enables us to deduce that $P \in \Ass_R H^2_J (M) \cup \Ass_R H^2_K (M)$. Consequently,
     $$
     \Ass_R H^2_I (M) \subset (\Ass_R H^2_I (M) \cap V(J+K)) \cup \Ass_R H^2_J (M) \cup \Ass_R H^2_K (M),
     $$ which completes the proof as the set on the right hand side is finite. 
     \end{proof}
    
We illustrate Theorem \ref{five dimensional part huneke} as follows.
    
    \begin{Example}\label{example dim 5}\rm Let $x_1,x_2,x_3,x_4,x_5, x_6$ be six indeterminates over a field $k$ with ${\rm char}\,k\neq 2$, and
 $$R=k[x_1,x_2,x_3,x_4,x_5, x_6]_{(x_1,x_2,x_3,x_4,x_5, x_6)}/(q), \quad q=\sum_{i=1}^6x_i^2.$$ Note $R$ is a (five-dimensional) factorial domain by \cite[Theorem 8.2]{Samuel}. Choose four $k$-linearly independent linear forms, which for simplicity we assume to be $x_1, x_2, x_3, x_4$. Consider the ideals $$I=(x_2 x_3, x_2x_4)R \quad \mbox{and} \quad J=(x_1^n, x_1 x_2, x_1 x_3)R, \quad n\geq 2.$$ Since $\Ht J=1$, we have $\dim R/J=4$. Furthermore, the image of $x_1^n$ in $R$ (hence an element of $J$) is $R/I$-regular. By Theorem \ref{five dimensional part huneke}, we conclude that $\Ass_R H^i_I (R/J)$ is finite for all $i$. 
 \end{Example}
    
\begin{rem}
    
  For more than one reason, Example \ref{example dim 5} cannot be deduced from \cite[Theorem 2.11]{Marley1}. Indeed,  $R$ is not regular and  $R/J$ is $R$-torsion.
\end{rem} 
    
\section{On the six-dimensional case}\label{sixdim}


This last section contains a theorem which we believe to be the first existing result towards the six-dimensional case of Huneke's conjecture. As expected, some conditions are needed. While most of the techniques used in the previous sections for the factorial case are employed here as well, we need to restrict ourselves to the smaller class of regular local rings containing a field. Such a constraint is essentially justified by the fact that the finiteness of $\Ass_R H^i_I (R)$ is well-known if $R$ belongs to this class of rings (a fact first proved in positive characteristic in \cite{HuS} and in characteristic zero in \cite{Lyu0}; see the lemma below), and the key strategy in our proof is to transfer somehow the finiteness of $\Ass H^i_I (R)$ to that of $\Ass_R H^i_I (M)$. 

\begin{Lemma}{\rm (\cite[Corollary 3.6(c)]{Lyu0})} \label{finiteass}
 If $R$ is a regular local ring containing a field, then $\Ass_R H^i_I (R)$ is finite for all $i$.   
\end{Lemma}

Below we establish two additional auxiliary lemmas.

\begin{Lemma}\label{lema tecnico dim 6}
 Let $P$ be a prime ideal of $R$ such that $R_P$ is analytically irreducible {\rm (}e.g., if $R_P$ is regular{\rm )}. If $\Ht P =n$ and $P \in \Supp_R H^n_I (R)$, then $P \in \Min_R R/I$.     
\end{Lemma}
\begin{proof} First, as a matter of standard notation, we denote the adic completion of a local ring $S$ with respect to its maximal ideal by $\widehat{S}$. Since formation of local cohomology commutes with localization, we have $H^n_{I_P} (R_P) \neq 0$. By hypothesis, the local ring $\widehat{R_P}$ is a domain. Therefore, the well-known Hartshorne-Lichtenbaum vanishing theorem (see \cite[Theorem 8.2.1]{Brodman}) yields
    $$
    \dim \widehat{R_P}/I_P\widehat{R_P} = 0,
    $$ or equivalently, $P_P\widehat{R_P}$ is a minimal prime over $I_P\widehat{R_P}$. Now, suppose the existence of a prime ideal $Q$ of $R$ satisfying
    $I \subset Q \subsetneqq P$.
    In particular,
    $I_P \subset Q_P \subsetneqq P_P$.
    Since $R_p \rightarrow \widehat{R_P}$ is faithfully flat, we can use \cite[Theorem 7.5(ii)]{Matsumura} to deduce that
    $$
    P_P\widehat{R_P} \cap R_P = P_P.
    $$ Furthermore, the flatness of the map $R_p \rightarrow \widehat{R_P}$ implies that it satisfies the going-down property (see \cite[Theorem 9.5]{Matsumura}). It follows that there is a prime ideal $Q' \subset P_P\widehat{R_P}$ such that 
    $Q' \cap R_P = Q_P$. Consequently,
    $$
    I_P\widehat{R_P} \subset Q_P\widehat{R_P} \subset Q^{'},
    $$
    which forces the equality $Q'=P_P\widehat{R_P}$.   
    Finally, we can write
    $$
    P_P=P_P\widehat{R_P} \cap R_P = Q' \cap R_P = Q_P,
    $$
    which is a contradiction. This shows that $P \in \Min_ R R/I$.
\end{proof}

Recall the ideal $I$ is {\it height-unmixed}
 if the associated prime divisors of $I$ have the same height. This holds for example if $R/I$ is a Cohen-Macaulay local ring.

\begin{Lemma}\label{lema tecnico 2}
Let $P$ be a prime ideal of $R$ with $\Ht P=n+1$ and assume that, for each prime ideal $Q\subset P$ with $\Ht Q \geq n$, the ring $R_Q$ is analytically irreducible {\rm (}e.g., if $R$ is locally regular{\rm )}. Suppose the minimal prime divisors of $I$ have the same height {\rm (}e.g., if $I$ is height-unmixed{\rm )} and $\Ht I < n$. If $P \in \Supp_R H^n_I (M)$ for some
finite $R$-module $M$, then  $P \in \Min_R H^n_I (R)$.
\end{Lemma}
\begin{proof}
Let us prove first that $P \in \Supp_R H^n_I (R)$. Suppose, by way of contradiction, that $H^n_{I_P} (R_P)=0$.
Since $P$ is not a minimal prime divisor of $I$, we must have $H^{n+1}_{I_P} (R_P)=0$ by Lemma \ref{lema tecnico dim 6}.
Now, we may pick suitable finite $R$-modules $N$, $N'$ fitting into short exact sequences
$$
0 \to N \to R^t \to M \to 0
\quad \mbox{and} \quad 0 \to N' \to R^s \to N \to 0$$
which in turn yield, respectively, the exact sequences in local cohomology 
\begin{equation}\label{exseq1}
H^n_{I_P} (R^t_P) \to H^n_{I_P} (M_P) \to H^{n+1}_{I_P} (N_P),    
\end{equation}
\begin{equation}\label{exseq2}
H^{n+1}_{I_P} (R^s_P) \to H^{n+1}_{I_P} (N_P) \to H^{n+2}_{I_P} (N'_P).    
\end{equation}
Since ${\rm dim}\,N'_P\leq {\rm dim}\,R_P=\Ht P=n+1$, we obtain $H^{n+2}_{I_P} (N'_P)=0$. Hence, it follows from \eqref{exseq2} that $H^{n+1}_{I_P} (N_P)=0$. Now, from \eqref{exseq1}, we get $H^n_{I_P} (M_P)=0$, which is a contradiction.

Finally, it suffices to argue that $P$ is a minimal element of $\Supp_R H^n_I (R)$. Suppose there exists an ideal $ Q \in \Supp_R H^n_I (R) $ with $Q \subsetneqq P$. Then, necessarily, $\Ht Q=n$. Now, it follows immediately from Lemma \ref{lema tecnico dim 6} that $Q \in \Min_R R/I$, which gives $\Ht Q= \Ht I<n$, a contradiction.
\end{proof}

For the last result of the paper, we use for a finite $R$-module $M$ and each $i\geq 0$  the notation $$\Ass_R^i M= \{P \in \Ass_R M \mid \Ht P =i\}.$$

\begin{Theorem}\label{theorem dim 6}
    Let $(R, {\bf m})$ be a six-dimensional regular local ring containing a field and $M$ a finite $R$-module whose annihilator is not contained in any minimal prime divisor of $I$ {\rm (}e.g., if ${\rm ann}_RM$ contains an $R/I$-regular element{\rm )}. In addition, suppose $\Ann_R M \not \subset {\bf m}^2$. If $I$ is height-unmixed, then $\Ass_R H^i_I (M)$ is finite for all $i \neq 2$.
\end{Theorem}
\begin{proof}
First, note  $\Ann_R M \neq 0$ and $R$ is, in particular, a domain. So, $\dim M=\dim R/{\rm ann_RM}\leq \dim R-1 = 5$. If $\dim M \leq 3$ then \cite[Corollary 2.7]{Marley1} yields the result. Now, suppose $\dim M = 4$. By \cite[Lemma 2.6]{Cuong}, the $R$-module $H^3_I (M)$ has finite support, hence finitely many associated primes. We can get the same conclusion concerning $H^4_I (M)$, by replacing $R$ with $R/{\rm ann_RM}$
and applying  \cite[Corollary 2.4]{Marley1}. Therefore, we may assume $\dim M =5$, in which case the same argument as above ensures  the finiteness of the support of $H^4_I (M)$ and $H^5_I (M)$.

So we need to show the finiteness of $\Ass_R H^3_I (M)$.
According to Lemma \ref{lemma dim menor que 2}, we may assume $\Ht I \leq 3$. Suppose first the case $\Ht I =3$. Given a prime ideal $Q \in \Supp_R H^i_I (M)$, we must have $\Ht Q \geq 3$. If $\Ht Q=3$, then $Q \in \Min_R R/I \cap V({\rm ann_RM})$, which gives a contradiction. Hence $\Ht Q \geq 4$, which yields $\dim H^i_I (M) \leq 2$ for all $i$, and once again we can apply Lemma \ref{lemma dim menor que 2}.

Notice that the case $\Ht I =1$ is easy to treat, because in this situation the ideal $I$ (being height-unmixed) must be principal, which puts us in a position to apply Lemma \ref{lemma principal ideal}.

Now we deal with the case $\Ht I = 2$, which is subtler.
Pick a prime ideal $P \in \Ass_R H^3_I (M)$. Then $\Ht P \geq 3$. If $\Ht P =3$, then [\citealp{Brodman}, Lemma 8.1.7] implies $P \in \Supp_R$ $H^3_I(R)$. Thus, as an application of Lemma \ref{lema tecnico dim 6}, we obtain that $P \in \Min_R R/I$, which is a contradiction.
So, $\Ht P >3$. At this point, we have shown the inclusion
$$
\Ass_R H^3_I (M) \subset \Ass_R^4 H^3_I (M) \cup \Ass_R^5 H^3_I (M) \cup \{{\bf m}\}.
$$
Moreover, since $I$ is a height-unmixed ideal with $\Ht I=2$, we may use Lemma \ref{lema tecnico 2} to obtain an inclusion $\Ass_R^4 H^3_I (M) \subset \Ass_R H^3_I (R)$, where the latter is finite by Lemma \ref{finiteass}. So, $\Ass_R^4 H^3_I (M)$ is finite as well.


It thus remains to prove the finiteness of the set $\Ass_R^5 H^3_I (M)$. With this in mind, we begin by claiming that $\Ann_R M   \setminus {{\bf m}}^2$ contains some $R/I$-regular element.
Indeed, otherwise we would have $$\Ann_R M \subset {{\bf m}}^2 \cup \displaystyle\bigcup_{Q \in \Ass_R R/I} Q $$ and then, by prime avoidance, either $\Ann_R M \subset {{\bf m}}^2$ or $\Ann_R M \subset Q$ for some $Q \in \Ass_R R/I$. The first possibility is ruled out by hypothesis, and so the second situation  holds. But notice that, as $I$ is height-unmixed, we have $\Ass_R R/I = \Min_R R/I$ and consequently $Q \in \Min_R R/I \cap \Supp_R M=\Min_R R/I \cap V({\rm ann_RM})$, which is a contradiction. This shows the claim.

Now, we can choose an $R/I$-regular element $$x \in \Ann_R M \setminus {{\bf m}}^2.$$ Set $T=R/(x)$, which is a five-dimensional regular local ring containing a field, and note $M$ has a natural structure of a finite $T$-module.
Moreover, by independence of the base ring (see \cite[Theorem 4.2.1]{Brodman}), the modules $H^3_I (M)$ and $H^3_{IT} (M)$ are isomorphic as $T$-modules.

Consider $P \in \Ass_R^5 H^3_I (M)$.
Then $PT \in \Ass_T H^3_{IT} (M)$. Furthermore, since $x \in \Ann_R M$ and $P \in \Supp_R M$, we get $x \in P$ so that $PT=P/(x)$ and consequently $$\Ht PT= \dim T - \dim T/PT =\dim T - \dim R/P =5-1=4.
$$ In addition, as $x$ is $R/I$-regular, we  deduce that $\dim R/(I+(x)) =3$ and so $\Ht IT=2$.

Let us assume first that the ideal $IT$ is height-unmixed.
In this case, we are in a position to apply Lemma \ref{lema tecnico 2}  which gives $PT \in \Ass_T H^3_{IT} (T)$, and this set must be finite by Lemma \ref{finiteass}. Hence, $\Ass_R^5 H^3_I (M)$ is also finite. Now, to deal with the general case, we can proceed in the same fashion as in the proof of Theorem \ref{four dimensional case} and write $IT=J \cap K$, where $J$ and $K$ are ideals of $T$ such that $J$ is height-unmixed with $\Ht J=2$, $\Ht K \geq 3$ and $\Ht (J+K) \geq 4$.
The Mayer-Vietoris sequence 
$$
H^3_{J+K} (M) \to H^3_J (M) \oplus H^3_K (M) \to H^3_{IT} (M) \to H^4_{J+K} (M)
$$
yields an inclusion
$$
\Ass^4_{T} H^3_{IT} (M) \subset \Supp_{T} H^3_{J+K} (M) \cup \Supp_{T} H^4_{J+K} (M) \cup \Ass^4_{T} H^3_J (M) \cup \Ass^4_{T} H^3_K (M).
$$
The finiteness of $\Ass^4_{T} H^3_J (M)$ is again a consequence of Lemma \ref{lema tecnico 2} together with Lemma \ref{finiteass}. The finiteness of $\Ass^4_{T} H^3_K (M)$ follows easily from Lemma \ref{lemma dim menor que 2} since $\dim T/K \leq 2$. Lastly,
$$
\Supp_{T} H^j_{J+K} (M) \subset \Min_T T/(J+K) \cup \{{\bf m}T\}
$$
which is therefore a finite set for each $j$.
This completes the proof of the finiteness of $\Ass^4_{T} H^3_{IT} (M)$, which implies that $\Ass^5_R H^3_I (M)$ is finite, as needed.
\end{proof}

\begin{rem}\rm An inspection of the above proof reveals that, unfortunately, our methods do not apply to the case ${\rm dim}\,R\geq 7$ (even when $R$ is a regular local ring containing a field). For instance, the (necessary) finiteness of some sets of primes  -- e.g., certain supports -- is no longer guaranteed, and in addition the use of Lemma \ref{lema tecnico 2} has a substantially weaker effect as it does not overcome the "gaps" that now appear in the argument due to the increase in ${\rm dim}\,R$.
\end{rem} 

\begin{rem} Since $H^i_I(M)=H^i_{\sqrt{I}}(M)$ for every $R$-module $M$ and $i\geq 0$, there is no loss of generality in assuming that $I$ is radical (thus yielding, in particular, the equivalence of the hypotheses "${\rm ann}_RM$ contains an $R/I$-regular element'' and "${\rm ann}_RM$ is not contained in any minimal prime divisor of $I$'', which also appeared in the main results of Section \ref{twothm}).
\end{rem}

The only case not touched by the last theorem leads us naturally to the following question.

\begin{quest}\label{2finite} \rm  Under the  hypotheses of Theorem \ref{theorem dim 6}, is it true that $\Ass_R H^2_I (M)$ is finite?
\end{quest}

\begin{rem}\rm In regard to Question \ref{2finite} above, it has to be pointed out that in general it is possible for $\Ass_R H^2_I (M)$ to be infinite. One explicit example, in the case where $\Ann_R M \subset {\bf m}^2$ (thus violating one of the key hypotheses of our Theorem \ref{theorem dim 6}), can be seen by applying the independence property of the base ring to \cite[Corollary 1.3]{Moty}.
\end{rem} 

We close the paper with an illustration of Theorem \ref{theorem dim 6}.

\begin{Example}\label{example dim 6}\rm  Let $x_1,x_2,x_3,x_4,x_5, x_6$ be six formal indeterminates over a field $k$, and consider the regular local ring $R=k[\![x_1,x_2,x_3,x_4,x_5,x_6]\!]$. Consider the $R$-ideal $$I=(x_1,x_2) \cap (x_3,x_4)$$ as well as the $R$-module $M=R/J$, where $J=(x_1 x_2, x_5)$. Notice that the indeterminate $x_5$ lies in $J\setminus {\bf m}^2$ and is $R/I$-regular. Thus, we can apply Theorem \ref{theorem dim 6} in order to deduce the finiteness of $\Ass_R H^i_I (R/J)$ for all $i \neq 2$.
\end{Example}
    
\bigskip

\noindent{\bf Acknowledgements.} The second-named author was partially supported by CNPq (grants 406377/2021-9 and 313357/2023-4).




\begin{thebibliography}{10}

     
    
\bibitem{Bah2} K. Bahmanpour, R. Naghipour, {\it Cofiniteness of local cohomology modules for ideals of small dimension}, J.
Algebra {\bf 321} (2009), 1997-2011.
    
    \bibitem{multi} B. Bhatt, M. Blickle, G. Lyubeznik, A. K. Singh, W. Zhang, {\it Local cohomology modules of a smooth ${\mathbb Z}$-algebra
have finitely many associated primes}, Invent. Math. {\bf 197} (2014), 509-519.
    

    \bibitem{Brodman} M. P. Brodmann, R. Y. Sharp,  {\it Local Cohomology. An Algebraic Introduction with Geometric Applications}, 2nd ed., Cambridge Stud. Adv. Math. {\bf 136}, Cambridge Univ. Press, Cambridge, 2013. 
    

    
    \bibitem{Goto} N. T. Cuong, S. Goto, N. V. Hoang, {\it On the cofiniteness of generalized local cohomology modules},  Kyoto J. Math. {\bf 55} (2015), 169–185.
    
    \bibitem{Cuong} N. T. Cuong, N. V. Hoang, {\it On the vanishing and the finiteness of supports of generalized local cohomology modules}, Manuscripta Math. {\bf 126} (2008), 59–72.

    \bibitem{Mafi1} K. Divaani-Aazar, A. Mafi,  {\it Associated primes of local cohomology modules of weakly Laskerian modules}, Comm. Algebra {\bf 34} (2006), 681–690.
    

\bibitem{Hu} C. Huneke, {\it Problems on local cohomology}, in: Free resolutions in commutative algebra and algebraic geometry, Sundance, Utah, 1990, Res. Notes Math. {\bf 2}, Jones and Bartlett, Boston, MA, 1992, pp. 93-108.


\bibitem{Hu-Koh} C. Huneke, J. Koh, {\it Cofiniteness and vanishing of local cohomology modules}, Math. Proc. Cambridge
Philos. Soc. {\bf 110} (1991), 421-429.

\bibitem{HuS} C. Huneke, R. Y. Sharp, {\it Bass numbers of local cohomology modules}, Trans. Amer. Math. Soc. {\bf 339} (1993), 765-779.

\bibitem{Moty} M. Katzman, {\it An example of an infinite set of associated primes of a local cohomology module}, J. Algebra
{\bf 252} (2002), 161-166.

    
    \bibitem{Kawazaki} K. I. Kawazaki, {\it Cofiniteness of local cohomology modules for principal ideals}, Bull. London Math. Soc. {\bf 30} (1998), 241-246

\bibitem{Lyu0} G. Lyubeznik, {\it Finiteness properties of local cohomology modules $($an application of D-modules to commutative algebra$)$}, Invent. Math. {\bf 113} (1993), 41-55.

\bibitem{Lyub} G. Lyubeznik, {\it A partial survey of local cohomology}, in: Local Cohomology and Its Applications, Marcel Dekker, 2001.
   
   \bibitem{Marley1} T. Marley, {\it The associated primes of local cohomology modules over rings of small dimension}, Manuscripta Math. {\bf 104} (2001), 519–525.
   
   \bibitem{Marley} T. Marley, J. C. Vassilev, {\it Cofiniteness and associated primes of local cohomology modules},  J. Algebra {\bf 256} (2002), 180–193.

\bibitem{Matsumura} H. Matsumura, {\it Commutative Ring Theory}, Cambridge Stud. Adv. Math. {\bf 8}, Cambridge Univ. Press, 1986.
    
   \bibitem{Mel} L. Melkersson, {\it Modules cofinite with respect to an ideal},  J. Algebra {\bf 285} (2005), 649–668.
    
    
   \bibitem{Roshan} H. Roshan-Shekalgourabi, M.Hatamkhani, {\it Hartshorne's question and weakly cofiniteness}, Math. Rep. (Bucur.) {\bf 22(72)} (2020), 329–340.
    
    
   \bibitem{Samuel} P. Samuel, {\it Lectures on Unique Factorization Domains}, Tata Institute of Fundamental Research, Bombay, 1964.
    
            
   \bibitem{Saremi} H. Saremi, A. Mafi, {\it On the Finiteness Dimension of Local Cohomology Modules}, Algebra Colloq. {\bf 21} (2014), 517–520.
    
   \bibitem{Bah3} M. Sedghi, K. Bahmanpour, R. Naghipour, {\it On the finiteness properties of local cohomology modules for regular local rings}, Tokyo J. Math. {\bf 40} (2017), 83–96.

\bibitem{Sin} A. K. Singh, {\it $p$-torsion elements in local cohomology modules}, Math. Res. Lett. {\bf 7} (2000), 165-176.
    
   \bibitem{SS} A. K. Singh, I. Swanson, {\it Associated primes of local cohomology modules and of Frobenius powers}, Int. Math. Res. Not. {\bf 33} (2004), 1703-1733.
   
    

\end{thebibliography}
\end{document}